%% Plain Tex file

%FINAL COPY OF ASKEY PAPER !!!!!!!!!!
\magnification=1200
\baselineskip=12truept
\vsize=8.5truein
%\vsize=10.125truein
\hsize=6.125truein
\overfullrule0pt
%\nopagenumbers

\def\nhang{\hangindent=4pc\hangafter=1}

\def\sqr#1#2{{\vcenter{\vbox{\hrule height.#2pt
        \hbox{\vrule width.#2pt height#1pt \kern#1pt
           \vrule width.#2pt}
            \hrule height.#2pt}}}}

\def\dotting{\leaders\hbox to 1em{\hfil.\hfil}\hfil}

\vglue.5truein
\noindent {\bf AN ANALOG OF THE FOURIER TRANSFORMATION}

\noindent{\bf FOR A ${\bf q}$-HARMONIC OSCILLATOR}
\vglue.65truein

\hskip 1.0truein {R.\ Askey,$^{\rm 1}$ N.\ M.\ Atakishiyev,$^{\rm
2}$ and
S.\ K.\ Suslov$^{\rm 3}$}

\bigskip
\hskip 1.0truein{$^{\rm 1}$Dept. of Math., U. of Wis.,
Madison, WI 53706, USA}

\hskip 1.0truein{$^{\rm 2}$Physics Institute, Baku 370143,
Azerbaijan}

\hskip 1.0truein{$^{\rm 3}$Kurchatov Institute of Atomic Energy,
Moscow
123182, Russia}
\vglue.65truein
\noindent {\bf ABSTRACT}
\medskip
A $q$-version of the Fourier transformation and some of its
properties
are discussed.

\bigskip

\noindent {\bf INTRODUCTION}
\medskip

Models of $q$-harmonic oscillators are being developed in
connection with quantum
groups and their various applications (see, for example, Refs.\
[M], [Bi], [AS1],
and [AS2]).  For a complete correspondence with the
quantum-mechanical
oscillator problem, these models need an analog of the Fourier
transformation
that relates the coordinate and momentum spaces.  In the present
work we
fill this gap for one of the models, the one based on the
continuous
$q$-Hermite polynomials [M], [AS1] when $-1<q<1$.

In Section I we assemble all those formulas from [AS1], which are
necessary for
the subsequent exposition.  In Section II we discuss the relation
between
the Mehler bilinear generating function for Hermite polynomials and
the
kernel $\exp \left({i\over h} px\right)$ of the Fourier
transformation
that connects the coordinate $x$ and momentum $p$ spaces [W].  We
used the
bilinear formula of L.\ J.\ Rogers to obtain a reproducing kernel
and an analogue of
the Fourier transform in the setting determined by the continuous
$q$-Hermite
polynomials of Rogers. Some properties of the $q$-Fourier
transformation
are discussed in Sections III-IV.

In the following we take $0<q<1$, although most of the formulas
remain correct
when $-1<q<0$. The limiting case $q\to 0$ is of some mathematical
interest.
\medskip

{\bf 1. ${\bf q}$-Hermite functions}.  The continuous $q$-Hermite
polynomials
were introduced by Rogers [R].  They can be defined by the three
term
recurrence relation
$$
2xH_n(x\mid q)=H_{n+1}(x\mid q)+(1-q^n)H_{n-1}(x\mid q) ,
\eqno(1.1)
$$
$H_0(x\mid q)=1$, $H_1(x\mid q)=2x$.  They are orthogonal
on $-1\le x=\cos\theta\le 1$ with
respect to a positive measure $\rho (x)$

$$\eqalignno{
\rho (x) &=4\sin\theta (qe^{2i\theta},
      qe^{-2i\theta};q)_\infty&(1.2)\cr
  &=4\sqrt{1-x^2}\prod^\infty_{j=1}[1-2(2x^2-1)q^j+q^{2j}]\cr}
$$
where the usual notations (see [GR]) are
$$
\eqalignno{
(a;q)_\infty& =\prod^\infty_{n=0} (1-aq^n)\, ,\cr
&&(1.3)\cr
(a,b;q)_\infty&=(a;q)_\infty\,(b;q)_\infty\,.\cr}
$$

A $q$-wave function for this $q$-harmonic oscillator was introduced
explicitly in [M], [AS1] as
$$
\psi_n(x)=\alpha d_n^{-1}\rho^{1/2}(x) H_n(x\mid q),\quad\alpha
=\left(
{1-q\over 2}\right)^{1/4}\, ,\eqno(1.4)
$$
where $d_n>0$ and
$$
d^{-2}_n={1\over 2\pi}(q^{n+1};q)_\infty\, .\eqno(1.5)
$$
These $q$-wave functions (1.4) satisfy
$$
\int^1_{-1}\psi_n(x)\psi_m(x)\,dx=\alpha^2\delta_{nm}\, .\eqno(1.6)
$$
See [AI] for a proof of this orthogonality relation.

$q$-annihilation $b$ and $q$-creation $b^+$ operators were
introduced
explicitly in
[M] and [AS1].  They satisfy the commutation rule
$$
bb^+-q^{-1}b^+b=1\eqno(1.7)
$$
and act on the $q$-wave functions defined in (1.4) by
$$
b\psi_n(x)=\tilde e_n^{1/2}\psi_{n-1}(x),\quad
b^+\psi_n(x)=\tilde e_{n+1}^{1/2}\psi_{n+1}(x)\eqno(1.8)
$$
where
$$
\tilde e_n={1-q^{-n}\over 1-q^{-1}}=q^{1-n}e_n\, .\eqno(1.9)
$$

In the limit when $q\to 1^-$ the functions $\psi_n(\alpha^2\xi )$
converge
to the classical wave functions
$$
\Psi_n(\xi )={1\over \pi^{1/4}(2^nn!)^{1/2}} H_n(\xi
)e^{-\xi^2/2}\eqno(1.10)
$$
 of the linear harmonic oscillator in the coordinate
representation.
\medskip

{\bf 2. An analog of the Fourier transformation}.  In proving that
the
Hermite functions (1.10) are complete in the space $L_2$ over
$(-\infty ,
\infty )$ and that a Fourier transform of any function from $L_2$
belongs to the
same space, the important role is played by the bilinear generating
function
(or the Poisson kernel) [W]
$$\eqalignno{
K_t(\xi ,\eta )&=\sum^\infty_{n=0} t^n\Psi_n(\xi )\Psi_n(\eta )\cr
  &=[\pi (1-t^2)]^{-1/2}\exp \left[{4\xi\eta t-
  (\xi^2+\eta^2)(1+t^2)\over 2(1-t^2)}\right]\, .&(2.1)\cr}
$$
Since in the limit $\gamma\to\infty$ the continuously
differentiable
function
$$
\delta (x,\gamma )={\gamma\over\sqrt\pi} e^{-\gamma^2x^2}
$$
becomes the Dirac delta function $\delta (x)$, from (2.1) in the
limit
$t\to 1^-$ follows the completeness property of the system (1.10),
$$
\sum^\infty_{n=0}\Psi_n(\xi )\Psi_n(\eta )=\delta (\xi -\eta )\,
.\eqno(2.2)
$$

On another hand, setting $t=i$ reduces the right-hand side of (2.1)
to the
kernel of the Fourier transformation, i.e.
$$
K_i(\xi ,\eta )={1\over\sqrt{2\pi}}e^{i\xi\eta}\, .
$$
Actually this idea provides the possibility of finding the Fourier
transformation
for difference analogs of the harmonic oscillator and its
$q$-generalizations,
when a priori it is not clear how one can define the explicit form
of the
kernels and how these transformations can look like.
As first examples  we considered  Kravchuk and Charlier functions
[AAS].  Here we shall give the necessary formulas for the
$q$-Hermite
functions (1.4).

Let us define a $q$-analog of (2.1) as
$$
K_t(x,y;q)=\alpha^{-2}\sum^\infty_{n=0} t^n\psi_n(x)\psi_n(y)\,
.\eqno(2.3)
$$
Substituting the formulas (1.4) and (1.5) into (2.3), we get
$$
K_t(x,y;q)={(q;q)_\infty\over 2\pi}\rho^{1/2}(x)\rho^{1/2}(y)
\sum^\infty_{n=0}{t^n\over (q;q)_n}H_n(x|q)H_n(y|q).\eqno(2.4)
$$
This series (the Poisson kernel) was summed as an infinite  product
by
Rogers [R], and the value is stated in [AI].  A very simple proof
was
given by Bressoud [Br].  The formula is
$$\eqalignno{
&\sum^\infty_{n=0} {t^n\over (q;q)_n} H_n(\cos\theta
|q)H_n(\cos\varphi |q)=\cr
&\qquad =(t^2;q)_\infty (te^{i(\theta +\varphi )}, te^{i(\theta
-\varphi)},
te^{-i(\theta +\varphi )}, te^{-i(\theta -\varphi
)};q)^{-1}_\infty\, ,&(2.5)\cr}
$$
so we can write (2.4) as
$$
K_t(x,y;q)={(q,t^2;q)_\infty \rho^{1/2}(x)\rho^{1/2}(y)\over
2\pi (te^{i(\theta +\varphi )}, te^{i(\theta -\varphi )},
te^{-i(\theta +\varphi )},
te^{-i(\theta -\varphi )};q)_\infty}\, .\eqno(2.6)
$$

In complete analogy with the case of the harmonic oscillator,
$$
\lim_{t\to 1^-}
K_t(x,y;q)=\alpha^{-2}\sum^\infty_{n=0}\psi_n(x)\psi_n(y)
=\delta (x-y)\, .\eqno(2.7)
$$
This is even easier to prove than in the classical case, since any
set of
polynomials orthogonal on a finite interval is complete in $L^2$,
and this is
equivalent to (2.7).  Formula (2.7) follows from
$$
t^n\psi_n(x)=\int^1_{-1}K_t(x,y;q)\psi_n(y)\,dy,\quad |t|<1\,
,\eqno(2.8)
$$
which is obvious by integration, which is justified by uniform
convergence.
An easy corollary of (2.3) is
$$
\int^1_{-1} K_t(x,y;q) K_\tau (y,x';q)\,dy=K_{t\tau}(x,x';q),\;\;
|t|<1, \;\; |\tau |<1\, .\eqno(2.9)
$$

Not only is the limit $t\to 1^-$ interesting, $t\to i$ is also
interesting.
The kernel in this case is the special case of (2.6) when $t=i$,
i.e.
$$
K_i(x,y;q)={4\over\pi}(q^2;q^2)_\infty
{[\sin\theta\sin \varphi (qe^{2i\theta}, qe^{-2i\theta},
qe^{2i\varphi},
qe^{-2i\varphi};q)_\infty ]^{1/2}\over
(ie^{i(\theta +\varphi )}, ie^{-i(\theta +\varphi )}, ie^{i(\theta
-\varphi )},
ie^{-i(\theta -\varphi )};q)_\infty}\,,\eqno(2.10)
$$
$x=\cos\theta$, $y=\cos\varphi$.

The classical Fourier transform has a kernel which is bounded, and
the
fourth power of it is the identity.  The fourth power of the
$q$-Fourier
transform,
$$
F_q [ \psi ] (x)=\lim_{r\to1^-}\int^1_{-1} K_{ir}(x,y;q)\psi
(y)\,dy
\,,\eqno(2.11)
$$
is also the identity, from (2.8) with $t\to i$ and $i^4=1$.  Since
$[-1,1]$
is compact, the kernel $K_i(x,y;q)$ can not be bounded and have
$F^4_q=I$.  It is not, and the singularity comes from the first
term
in the four products in the denominator of (2.10).  The singular
part
is just the value of (2.10) when $q=0$, i.e.
$$
{-(\sin\theta \sin\varphi )^{1/2}\over \pi\cos (\theta +\varphi
)\cos
(\theta -\varphi )}\, .\eqno(2.12)
$$
Thus this $q$-Fourier transform looks more like a weighted Hilbert
transform
than it looks like the classical Fourier transform. An explicit
form of the
transformation (2.11) is
$$\eqalignno{
&F_q[\psi](x)=Pv\int^1_{-1}K_i(x,y;q)\psi(y)\,dy&(2.13)\cr
&\qquad +{i\over2}\left[ k(x)\psi (\sqrt{1-x^2})+
k(-x)\psi(-\sqrt{1-x^2})\right]\,,\cr}
$$
where
$$
k(x)={x\over\left(x^2(1-x^2)\right)^{1/4}}\left( \prod^\infty
_{k=1} {1+4ix\sqrt{1-x^2}q^k-q^{2k}\over
1-4ix\sqrt{1-x^2}q^k-q^{2k}} \right)^{1/2}
$$
and $Pv$ denotes Cauchy's principal value integral.

The inverse of this transformation follows from
$$
\lim_{r,r'\to1^-}\int^1_{-1} K_{ir}(x,y;q) K^*_{ir'}(y,x';q)\,dy
=\delta (x-x')\, ,\eqno(2.14)
$$
where $*$ denotes the complex conjugate.  Formula (2.14) follows
from (2.9) when we observe that
$$
K^*_{ir}(y,x;q) =K_{-ir}(y,x;q)\, .\eqno(2.15)
$$
Another form of (2.14) is
$$
\lim_{r,r'\to 1^-}\int^1_{-1}dy\,K^*_{ir'}(x,y;q)\,
\int^1_{-1}K_{ir}(y,z;q)f(z)\,dz=f(x)\,.\eqno(2.16)
$$
Indeed, by (2.9) and (2.15) we get
$$\eqalign{
&\int^1_{-1}dy\,K^*_{ir'}(x,y;q)\,
\int^1_{-1}K_{ir}(y,z;q)f(z)\,dz\cr
&\qquad =\int^1_{-1}K_{rr'}(x,z;q)f(z)\,dz
\longrightarrow f(x)\cr }
$$
with $r,r'\to 1^-$.

By using the change of variables $x=\alpha^2\xi$ and
$y=\alpha^2\eta$, where $\alpha^2 = \left({1-q\over2}
\right)^{1/2}$, from (2.3) in the limit $q\to 1^-$ we can
write
$$\eqalign{
&\alpha^2K_t(\alpha^2\xi,\alpha^2\eta;q)
=\sum^\infty_{n=0}t^n\psi_n(\alpha^2\xi)
\psi_n(\alpha^2\eta) \cr
&\qquad\longrightarrow\sum^\infty_{n=0}
\Psi_n(\xi)\Psi_n(\eta)\,=K_t(\xi,\eta)\,,\cr}
$$
or
$$
\lim_{q\to 1^-}\alpha^2K_t(\alpha^2\xi,\alpha^2\eta;q)
=K_t(\xi,\eta)\,.\eqno(2.17)
$$
With the aid of the same consideration, from (2.11) we get
$$\eqalign{
F_q[\psi](x)=&\lim_{r\to 1^-}\int_{-1/\alpha^2}^{1/\alpha^2}\,
\alpha^2K_{ir}(\alpha^2\xi,\alpha^2\eta;q)\psi(\alpha^2\eta)\,d\eta\cr
&\longrightarrow{1\over\sqrt{2\pi}}\int^\infty_{-\infty}
e^{i\xi\eta}f(\eta)\,d\eta =F[f](\xi)\,,\cr}
$$
where $f(\eta) =\lim_{q\to 1^-}\psi(\alpha^2\eta)$. Therefore,
the classical Fourier transformation can be realized as a limiting
case of the $q$-Fourier transform.
\medskip

{\bf 3. ``Momentum and position'' operators}.  As is well known
from quantum
mechanics, the kernel of the classical Fourier transformation is
also an
eigenfunction of the momentum and position operators, i.e.
$$\eqalignno{
Q_\xi K_i(\xi ,\eta )&=\xi K_i(\xi ,\eta )\, ,\cr
&&(3.1)\cr
P_\xi K_i(\xi ,\eta )&=i^{-1} {d\over d\xi}\left(
{1\over\sqrt{2\pi}} e^{i\xi\eta}
\right) =\eta K_i(\xi ,\eta )\, .\cr}
$$
We can retain this property in the case of the $q$-Fourier
transformation,
if we choose the
corresponding operators for the kernel (2.10) as
$$\eqalignno{
Q &={\sqrt{1-q}\over 2} \left( q^{N/2}b+b^+q^{N/2}\right)\, ,\cr
&&(3.2)\cr
P &={\sqrt{1-q}\over 2i} \left( q^{N/2}b-b^+q^{N/2}\right)\, ,\cr}
$$
where $N=\log [1-(1-q^{-1})b^+b]/\log q^{-1}$ is the ``particle
number''
operator.
In fact, the action of the operator $Q$  on the kernel (2.10) gives
$$\eqalignno{
\alpha^2Q_xK_i(x,y;q)&={\sqrt{1-q}\over 2}
\left( q^{N_x/2}b_x+b_x^+q^{N_x/2}\right)
\sum^\infty_{n=0}i^n\psi_n(x)\psi_n(y)&(3.3)\cr
&={\sqrt{1-q}\over 2}\sum^\infty_{n=0}i^n\left[
e_n^{1/2}\psi_{n-1}(x)
   +e^{1/2}_{n+1}\psi_{n+1}(x)\right]\psi_n(y)\, .\cr}
$$
From the three-term recurrence relation (1.1) for the continuous
$q$-Hermite
polynomials, it follows that
$$
e_{n+1}^{1/2}\psi_{n+1}(x)+e_n^{1/2}\psi_{n-1}(x)={2\over\sqrt{1-
q}}\,x\,
\psi_n(x)\, .\eqno(3.4)
$$
Therefore for the ``position'' operator $Q$ we have, indeed,
$$
Q_xK_i(x,y;q)=xK_i(x,y;q)\, .\eqno(3.5)
$$

Exactly in the same manner, one can obtain directly from
definitions
(3.2) and expansion (2.3) for $t=i$, that
$$
P_xK_i(x,y;q)=Q_yK_i(x,y;q)\, ,\eqno(3.6)
$$
and, consequently,
$$
P_xK_i(x,y;q)=yK_i(x,y;q)\, .\eqno(3.7)
$$
Equations
$$\eqalignno{
&Q_xK_t(x,y;q) =xK_t(x,y;q)\,,\cr
&&(3.8)\cr
&P_xK_t(x,y;q) ={2ty-(1+t^2)x \over i(1-t^2)}
K_t(x,y;q)\cr}
$$
are an extension of (3.5) and (3.7).

The commutation rule of the operators (3.2) is
$$
[Q,P] = i\,{1 -q\over 2}\, q^N \eqno (3.9)
$$
and, therefore, the Hamiltonian of the $q$-oscillator,
$$
\tilde H = b^+ b = {1 - q^{-N}\over 1 - q^{-1} }\, , \eqno (3.10)
$$
has the following form
$$
\tilde H (P,Q) = {i - 2(1 - q)^{-1} [Q,P]\over 2q^{-1}
[Q,P]} \eqno (3.11)
$$
in terms of these momentum and position operators.

The equation
$$
P_x^2K_i(x,p;q) =p^2K_i(x,p;q)\eqno(3.12)
$$
may be considered as an equation of motion for a``$q$-free
particle".
\medskip

{\bf 4. Some properties of the $q$-Fourier transformation}.  The
kernel $K_i(\xi ,\eta )$ corresponds to the classical Fourier
transformation
$$
f(\xi )={1\over\sqrt{2\pi}}\int^\infty_{-\infty} e^{i\xi \eta}
g(\eta )\,d\eta =F[g](\xi )\, ,
\eqno(4.1)
$$
which has well-known properties [W].  With the aid of the kernel
(2.6) we
defined a $q$-version of the Fourier transformation by
$$
\psi (x)=\lim_{r\to 1^-}\int^1_{-1}K_{ir}(x,y)\varphi (y)\,dy
=F_q[\varphi ](x)\, .\eqno(4.2)
$$
We can establish the simple properties of this generalization.  The
orthogonality property (2.14) of the kernel (2.6) results in the
inversion
formula
$$
\varphi (y)=\lim_{r'\to 1^-}\int^1_{-1}K^*_{ir'}(x,y)\psi (x)\,dx\,
,\eqno(4.3)
$$
as well as in the relation $\| \varphi\|^2=\|\psi\|^2$.

Analogs of the properties
$$
\eqalignno{
&i^{-1} {d\over d\xi }F[g](\xi )=F[\eta g](\xi )\, ,\cr
&&(4.4)\cr
&iF\left[ {dg\over d\eta}\right](\xi )=\xi F[g](\xi )\, ,\cr}
$$
have the forms
$$
\eqalignno{
&P_xF_q[\varphi ](x)=F_q[y\varphi ](x)\, ,\cr
&&(4.5)\cr
&F_q[P_y\varphi ](x)=-xF_q[\varphi ](x)\, .\cr}
$$
Moreover, the following properties of the Fourier transform,
$$
\eqalignno{
&F[f](\xi+\xi_0)=e^{i\xi_0P_\xi}F[f](\xi)=F[e^{i\xi_0\eta}f](\xi)\,
,\cr
&&(4.6)\cr
&F[f(\eta-\eta_0)](\xi)= F[e^{-i\eta_0P_\eta}f](\xi)=
e^{i\xi\eta_0}F[f](\xi)\, ,\cr}
$$
admit a generalization
$$
\eqalignno{
&K_i(x_0,P_x)F_q[\varphi](x)=F_q[K_i(x_0,y)\varphi](x)\, ,\cr
&&(4.7)\cr
&F_q[K_i(y_0,-P_y)\varphi](x)=K_i(x,y_0) F_q[\varphi](x)\,.\cr}
$$

To define a $q$-version of the convolution $\, \varphi * \psi  \,
$ let
us retain the property
$$
F_q[ \phi * \psi ] = F_q  [\phi ] \cdot F_q[\psi ]\, ,\eqno
(4.8)
$$
or
$$
\eqalignno {
   &\lim_{r\to 1^-} \int^1_{-1} K_{ir} (x,z) (\varphi * \psi )(z)\,
dz = \cr
   &&(4.9)\cr
   &\qquad\lim_{r',r''\to 1^-} \int^1_{-1} \int^1_{-1}
   K_{ir'} (x,y)K_{ir''} (x,y')\varphi (y)
   \psi (y')\,dy\,dy'\, . \cr}
$$
Using (2.14) and (4.9) we arrive at the definition
$$
\eqalignno {
(\varphi * \psi )(z) & =\lim_{r,r',r''\to 1^-} \int^1_{-1}
\int^1_{-1}
dy\,dy'\, \varphi (y)\psi(y')\cr
&&(4.10)\cr
& \int^1_{-1} dx\,K_{ir'} (x,y)K_{ir''} (x,y')
                  K^*_{ir} (x,z)\, . \cr}
$$
The usual properties,
$$
\eqalign {
\varphi * \psi & = \psi * \varphi\, ,\cr
(\varphi * \psi ) * \chi & = \varphi * (\psi * \chi )\,,\cr}
$$
are valid.

In the limit $q \to 1^-$ we can write
$$
\eqalignno {
&\alpha^4 \int^1_{-1} K_{ir'} (x,y)K_{ir''} (x,y')
K^*_{ir} (x,z)\,dx \cr
&&(4.11)\cr
&\qquad\longrightarrow {1\over (2\pi)^{3/2}}
 \int^\infty_{-\infty} e^{i\xi (\eta + \eta '-\zeta )}\,
d\xi = {1\over \sqrt{2\pi}} \delta (\eta + \eta ' - \zeta )\,, \cr}
$$
and, therefore,
$$
(\varphi * \psi )(z) \longrightarrow {1\over\sqrt{2\pi}}
\int^\infty_{-\infty } f(\eta) g(\zeta - \eta )\,d\eta\, .
$$

We think the transformation (4.2) deserves a more detailed
consideration.
\bigskip
\bigskip

\baselineskip=12pt
\noindent{\bf REFERENCES}
\medskip

\nhang{[AAS] R.\ Askey, N.\ M.\ Atakishiyev and S.\ K.\ Suslov --
Fourier transformations
for difference analogs of the harmonic oscillator.  Proceedings of
the XV Workshop on
High Energy Physics and Field Theory, Protvino, Russia, 6--10 July
1992.}
\medskip
\nhang{[AI] R.\ Askey and M.\ E.\ H.\ Ismail -- A generalization of
the
ultraspherical polynomials.  Studies in Pure Mathematics (P.\
Erd\"os, ed.),
Birkh\"auser, Boston, Massachusetts, pp.\ 55--78, 1983.}
\medskip
\nhang{[AS1] N.\ M.\ Atakishiyev and S.\ K.\ Suslov -- Difference
analogs of
the harmonic oscillator.  Theoretical and Mathematical Physics,
Vol.\ 85,
No.\ 1, pp.\ 1055--1062, 1990.}
\medskip
\nhang{[AS2] N.\ M.\ Atakishiyev and S.\ K.\ Suslov -- A
realization of the
$q$-harmonic oscillator.  Theoretical and Mathematical Physics,
Vol.\ 87, No.\ 1,
pp.\ 442--444, 1991.}
\medskip
\nhang{[Bi] L.\ C.\ Biedenharn -- The quantum group
$\hbox{SU}_q(2)$ and a
$q$-analogue of the boson operators.  J.\ Phys.\ A: Math.\ Gen.,
Vol.\ 22,
No.\ 18, pp.\ L873--L878, 1989.}
\medskip
\nhang{[Br] D.\ M.\ Bressoud -- A simple proof of Mehler's formula
for $q$-Hermite
polynomials, Indiana Univ.\ Math.\ J., Vol.\ 29, pp.\ 577--580,
1980. }
\medskip
\nhang{[GR] G.\ Gasper and M.\ Rahman -- Basic Hypergeometric
Series.
Cambridge University Press, Cambridge, 1990.}
\medskip
\nhang{[M] A.\ J.\ Macfarlane -- On $q$-analogues of the quantum
harmonic
oscillator and the quantum group $\hbox{SU}(2)_q$.  J.\ Phys.\ A:
Math.\ Gen.,
Vol.\ 22, No.\ 21, pp.\ 4581--4588, 1989.}
\medskip
\nhang{[R] L.\ J.\ Rogers -- Second memoir on the expansion of
certain infinite
products. Proc.\ London Math.\ Soc., Vol.\ 25, pp.\ 318--343,
1894.}
\medskip
\nhang{[W] N.\ Wiener -- The Fourier Integral and Certain of Its
Applications.
Cambridge University Press, Cambridge, 1933.}
\bigskip

%\noindent Dept.\ of Mathematics, U.\ of Wisconsin, Madison, WI
~53706, USA

%\noindent Physics Institute, Baku 370143, Azerbaijan

%\noindent Kurchatov Institute of Atomic Energy, Moscow 123182,
Russia

\bye